\documentclass{article}
\usepackage{latexsym,amssymb,amsmath}

\begin{document}

\title{On simple endotrivial modules}
\author{Geoffrey R. Robinson,\\
Institute of Mathematics,\\
University of Aberdeen,\\
Aberdeen AB24 3UE,\\
Scotland\\
g.r.robinson@abdn.ac.uk }
\date{}
\maketitle

\begin{abstract}
We show that when $G$ is a finite group which contains an elementary
Abelian subgroup of order $p^{2}$ and $k$ is an algebraically closed
field of characteristic $p,$ then the study of simple endotrivial
$kG$-modules which are not monomial may be reduced to the case
when $G$ is quasi-simple.

\end{abstract}

\section{Introduction:} Let $G$ be a finite group which contains an elementary Abelian $p$-subgroup of order $p^{2},$
and $k$ be an algebraically closed field of prime characteristic $p.$ The category of finite-dimensional $kG$-modules 
contains infinitely many isomorphism types of indecomposables, and this category is wild unless $p=2$ and $G$ has a Sylow 
$2$-subgroup which is dihedral or semidihedral.

\medskip
E.C. Dade introduced the notion of an endotrivial $kG$-module as an example of a type of $kG$-module which
occurs naturally in some representation theoretic problems, but may be more amenable to a reasonable
classification. A finite dimensional $kG$-module $M$ is said to be endotrivial if the trivial
module $k$ is a direct summand of $M \otimes M^{*}$ and is the the only non-projective summand of this
module. Endotrivial modules have been studied quite intensively in recent years (see, for example,
\cite{carlmazzthev} and the references therein). Notice that if $M$ is an endotrivial $kG$-module
with ${\rm dim}_{k}(M) >1,$ then no non-identity $p$-element of $G$ acts trivially on $M.$

\medskip
In Navarro-Robinson \cite{navrob}, a conjecture of J.F. Carlson, N. Mazza and J. Th\'evenaz
\cite{carlmazzthev} was proved, which completed the proof of the fact that if
$G$ as above is also $p$-solvable, then any simple endotrivial $kG$-module $V$ is $1$-dimensional. Here, we extend the
Clifford theoretic analysis of \cite{navrob} to reduce the study of simple endotrivial modules to the
case of groups which are close to quasi-simple with one type of exception in which the situation is very well-controlled.

\medskip
We will prove that if $G$ is an arbitrary finite group containing an
elementary Abelian $p$-subgroup of order $p^{2},$ and $V$ is a
faithful simple endotrivial $kG$-module, then there are just two
possibilities:

\medskip
\noindent i) $F^{*}(G) = Z(G)L,$ where $L$ is quasi-simple of order
divisible by $p$ and $L$ acts irreducibly on $V.$

\medskip
\noindent ii) $F^{*}(G) = O_{p^{\prime}}(G) = A$ is Abelian and there is a
strongly $p$-embedded subgroup $H$ of $G$ such that $V$ is induced from a $1$-dimensional $kH$-module.

\section{Some preliminaries and reductions}

\medskip
We first recall some now standard group-theoretic definitions and notation. A component of a finite group $G$
is a quasi-simple subnormal subgroup. Distinct components of $G$ centralize each other,
and the central product of all components of $G$ is denoted by $E(G).$ The generalized Fitting subgroup
of $G$ is denoted by $F^{*}(G)$ and is the (central) product $F(G)E(G).$ We always have
$C_{G}(F^{*}(G)) = Z(F(G)).$

\medskip
The proper subgroup $H$ of $G$ is said to be strongly $p$-embedded in $G$ if $p$ divides $|H|$
but $H \cap H^{g}$ is a $p^{\prime}$-group for all $g \in G \backslash H.$ Strongly
$p$-embedded subgroups play a large role in this paper, by virtue of Lemma 1, part iii),
of \cite{navrob}.

\medskip
\noindent {\bf Lemma 1:} \emph{Let $G$ be a finite group with a strongly $p$-embedded subgroup $H.$
Then $\cap_{g \in G} H^{g} \leq O_{p^{\prime}}(G),$ and equality holds
if $G$ contains a subgroup of type $(p,p).$ Furthermore, in the latter case,
$O_{p^{\prime},p}(G) = O_{p^{\prime}}(G).$}

\medskip
\noindent {\bf Proof:} Let $N = \cap_{g \in G}H^{g},$ let $P \in {\rm Syl}_{p}(H)$
and set $Q = P \cap N.$ If $ Q \neq 1,$ then $G = NN_{G}(Q) \leq H.$ Hence
$|N|$ is not divisible by $p.$ If $G$ contains a subgroup $R$ of type $(p,p),$
then $$O_{p^{\prime}}(G) \leq \langle C_{G}(x): x \in R^{\#} \rangle \leq H.$$

\medskip
In that case, set $N = O_{p^{\prime},p}(G).$ If $N > O_{p^{\prime}}(G),$
then $N = O_{p^{\prime}}(G)(P \cap N) \leq H,$ a contradiction as $N \lhd G$
and $p$ divides $|N|.$

\medskip
\noindent {\bf Remark:} If $G$ is a finite group with a strongly
$p$-embedded subgroup $H$ and $O_{p^{\prime}}(G) = 1$ then $F(G) =
1,$ and $G$ has a unique minimal normal subgroup, which is simple.
This is well-known, but notice that if $M$ is a minimal normal
subgroup of $G,$ but is not simple, then $G = MH$ by a Frattini
argument, while if $L,K$ are any two of the simple direct factors of
$M,$ then $L \leq C_{G}(K) \leq C_{G}(P \cap K) \leq H.$ Hence $M \leq H$ as $L$ was arbitrary, a
contradiction.

\medskip
\noindent {\bf Lemma 2:} \emph{Let $G$ be a finite group, $V$ be an
endotrivial simple $kG$-module, and let $N$ be a subnormal subgroup
of $G$ of order divisible by $p.$ Then ${\rm Res}^{G}_{N}(V)$ is
simple. In particular, if $N$ is a component of $G,$ and $V$ is
faithful, but not $1$-dimensional, then $C_{G}(N) = Z(G),$ so
$N$ is the unique component of $G.$ }

\medskip
\noindent {\bf Proof:} The first claim is is clear if ${\rm dim}_{k}(V) = 1,$ so
suppose that this is not the case. By induction, it suffices to
treat the case that $N$ is normal, so suppose that this is the case.
Then since $V$ is endotrivial, we may suppose that $G$ acts
faithfully on $V,$ since the kernel of the action of $G$ on $V$ is a
$p^{\prime}$-group. Let $U$ be a simple summand of ${\rm
Res}^{G}_{N}(V)$ and set $H =I_{G}(U) \geq N.$ If $H < G,$ then $V$
is induced from a simple endotrivial $kH$-module. Since $V$ is
endotrivial, it follows from Lemma 1 of Navarro-Robinson
\cite{navrob} that $H$ is strongly $p$-embedded, contrary to Lemma 1
above, as p divides $|N|.$  Hence $H = G$.

\medskip
Now there is a $p^{\prime}$-central extension ${\hat G}$ of $G$ such
that $U$ extends to a $k{\hat G}$-module and ${\hat G}$ has a normal
subgroup ${\hat N}$ which may be identified with $N,$ and $V \cong U
\otimes W,$ where $W \cong {\rm Hom}_{k{\hat N}}(U,V).$ Now both $U$
and $W$ are endotrivial by Lemma 1 of \cite{navrob}. Since $N$ acts
trivially on $W,$ some non-trivial $p$-element of ${\hat G}$ acts
trivially on $W,$ so that $W$ must be $1$-dimensional. Hence $V$
restricts irreducibly to $N,$ so $C_{G}(N) = Z(G).$

\medskip
We state the next Lemma in more generality than needed here, since it
may have application and interest beyond the question at hand.

\medskip
\noindent {\bf Lemma 3:} \emph{Let $V$ be a primitive indecomposable
endotrivial $kG$-module, where $G$ is a finite group which contains
an elementary Abelian subgroup $Q$ of order $p^{2}.$ Then
$O_{p^{\prime}}(G)$ acts as scalars on $V$.}

\medskip
\noindent {\bf Proof:} Let $U$ be a simple summand of ${\rm
Res}_{O_{p^{\prime}}(G)}^{G}(V).$ Since $V$ is primitive, $I_{G}(U)=G.$
Now there is a finite $p^{\prime}$-central extension
${\hat G}$ such that $U$ extends to a $k{\hat G}$-module, and as
such $V \cong U \otimes W,$ where $W = {\rm Hom}_{kO_{p^{\prime}}({\hat G})}(U,V)$
(via $u \otimes f \to uf$ for $ u \in U, f \in {\rm Hom}_{kO_{p^{\prime}}({\hat G})}(U,V)).$
 Then $V$ is still an endotrivial $k{\hat G}$-module, so $U$ and $W$ are both
endotrivial $k{\hat G}$-modules by Lemma 1 of \cite{navrob}. Since ${\hat G}$
contains a subgroup ${\hat Q}$ isomorphic to $Q$ and since $U$ is an endotrivial
simple $kO_{p^{\prime}}({\hat G}){\hat Q}$-module, $U$ must be
$1$-dimensional by the main result of \cite{navrob}.

\medskip
\noindent {\bf Lemma 4:} \emph{Let $V$ be a faithful imprimitive simple endotrivial $kG$-module,
where $G$ is a finite group containing an elementary Abelian subgroup of order $p^{2}.$
Then $A = O_{p^{\prime}}(G)$ is Abelian and either $A = C_{G}(A)$ or $A = Z(G).$ In the latter case,
$G$ has a component $L$ of order divisible by $p$ which acts irreducibly on $V.$}

\medskip
\noindent {\bf Proof:} $V$ is induced from a primitive simple
$kH$-module, say $U$  of some proper subgroup $H$ of
$G.$  By Lemma 1 of \cite{navrob}, $H$ is strongly $p$-embedded in $G.$
Then $A = O_{p^{\prime}}(G) \leq N = \cap_{g \in G}H^{g}.$ By
Lemma 3 (applied within $H,$ as $U$ is primitive),  $A$ acts as scalars on $U.$ Hence $A$ is Abelian by
Clifford's Theorem. If $C= C_{G}(A) >A,$ then $|C|$ is divisible by
$p,$ so that $C$ acts irreducibly on $V$ by Lemma 2. In that case, $A$ acts
as scalars on $V$ and $A \leq Z(G).$ Thus $G$ has a component $L$ of
order divisible by $p$ and Lemma 2 applies.

\medskip
Hence we may suppose that $C = A.$ Now, by Lemma 1 and the ensuing discussion, $G/A$ has a corefree strongly $p$-embedded subgroup,
and $F^{*}(G/A)$ is simple.

\section{Proof of the main Theorem}

\medskip
\noindent {\bf Theorem:} \emph{ Let $G$ be a finite group containing an elementary Abelian
subgroup of order $p^{2},$ and let $V$ be a faithful simple endotrivial $kG$-module. Then
one of the following occurs:}

\medskip
\noindent \emph{ i) $E = E(G)$ is quasi-simple and acts irreducibly on $V$.}

\medskip
\noindent \emph{ii) $O_{p^{\prime}}(G) = A$ is Abelian, $C_{G}(A) = A$ and $V$ is induced
from a $1$-dimensional module of a strongly $p$-embedded subgroup $H$ of $G.$}

\medskip
\noindent {\bf PROOF:} If $V$ is primitive, then $O_{p'}(G)$ acts as scalars on $V$ by Lemma 3. Hence
$G$ must have a component $L$ of order divisible by $p$ and case i) holds by Lemma 2.

\medskip
We may thus suppose that $V$ is imprimitive, and by Lemma 2, that $G$ has no component
of order divisible by $p.$ Then $V$ is induced from a primitive simple endotrivial
$kM$-module $W$ for some strongly $p$-embedded subgroup $M$ of $G.$ If $W$
is $1$-dimensional, we are done, so suppose that ${\rm dim}_{k}(W) > 1.$
Set $A = O_{p^{\prime}}(G).$ By Lemma 4, since $G$ has no component of order divisible
by $p,$ we must have $A = F^{*}(G)$ and $C_{G}(A) = A.$ Also $A \leq M.$

\medskip
Applying the earlier Lemmas to $M$ in its action on $W,$ we know that the kernel $K$ of
the action of $M$ on $W$ is a $p^{\prime}$-group, that $O_{p^{\prime}}(M)/K$ acts as scalars on $W,$
and that $M/K$ has a unique component which acts irreducibly on $W.$

\medskip
Set ${\bar G} = G/A,$ etc.. Then ${\bar M}$ is still strongly $p$-embedded in ${\bar G}.$
Also $A/K$ is central in $M/K,$ so that $M/K$ is a $p^{\prime}$-central extension of ${\bar M}.$
Since $M/K$ has a single component (of order divisible by $p$) and $O_{p}(M/K) = 1,$ we know that
${\bar M}$ also has a unique component and $O_{p}({\bar M}) =1.$

\medskip
The finite simple groups of order divisible by $p$ which contain an elementary Abelian subgroup of order $p^{2}$
and have a strongly p-embedded subgroup are carefully listed in 7.6.1 and 7.6.2 of (Gorenstein, Lyons and Solomon,
\cite{gorlysol}).  Since $O_{p}({\bar M}) =1$ most of these possibilities are excluded in the present situation.
Also, since ${\bar M}$ has a unique component of order divisible by $p,$
the possibility that $E({\bar G}) \cong A_{2p}$ is also excluded. The only possibility that remains is that $p =5,$
$E({\bar G}) \cong Fi_{22}$ and ${\bar M} \cap E({\bar G}) \cong {\rm Aut}(O_{8}^{+}(2)).$

\medskip
By Lemma 2, $W$ is a (faithful) simple endotrivial module of dimension greater than $1$ for a component, $L$ say,
of $M/K.$ This component $L$ is isomorphic to either $O_{8}^{+}(2)$ or a double cover of $O_{8}^{+}(2).$
From \cite{modatl} (also from \cite{gap}) no such group has a non-trivial absolutely irreducible module in characteristic $5$
whose dimension is congruent to $\pm 1$ (mod $25$). This contradicts the fact that $W$ is endotrivial of dimension greater
than $1,$ so that  ${\rm dim}_{k}(W)^{2} \equiv 1$ (mod $25$).

\medskip
\section{Concluding Remarks}
For the sake of completeness, we prove that simple monomial endotrivial $kG$-modules of dimension greater
than one exist for some groups $G$ which are not quasi-simple.

\medskip
\noindent {\bf LEMMA 5:} \emph{ Let $X$ be a finite group of order divisible by $p.$ 
Suppose that $O_{p^{\prime}}(X) = 1$ and that $X$ has a strongly $p$-embedded subgroup $Y.$ Then there is a finite group $G$ such that $F^{*}(G) = A$
is an Abelian $p^{\prime}$-group with $G/A \cong X$ and there is a simple endotrivial $kG$-module which is induced from a $1$-dimensional
module of a proper subgroup of $G.$}

\medskip
\noindent {\bf PROOF:} We might as well suppose that $Y$ is a maximal subgroup of $X,$ and we do so.
Let $q$ be a prime which does not divide $|X|$ and set $F = {\rm GF}(q).$ Then the permutation module 
${\rm Ind}_{Y}^{X}(F)$ has a non-trivial irreducible summand, say $A.$ By Frobenius reciprocity,
$Y$ fixes a non-identity element of the Abelian group $A,$ but $X$ fixes no such element.
By Brauer's permutation Lemma, $Y$ fixes a non-trivial linear character, say $\lambda,$ of $A,$
but $X$ fixes no such linear character of $A.$ Let $G$ be the semidirect product $AX.$ Then $H = AY$
is a strongly $p$-embedded maximal subgroup of $G$ and $H$ is the inertial subgroup $I_{G}(\lambda).$

\medskip
Let $U$ be an irreducible $kA$-module with Brauer character $\lambda.$ Now $A = [A,Y] \times C_{A}(Y),$
so that $U$ extends to a $1$-dimensional $kH$-module ${\tilde U}$ with $[A,Y]$ in its kernel.     
We still have $H = I_{G}(U).$ Let $V$ be a composition factor of $Ind_{H}^{G}({\tilde U}).$
Then $V$ covers the $kA$-module $U,$ so that $V$ is induced from a simple $kI_{G}(U)$-module.
Hence $[G:H]$ divides ${\rm dim}_{k}(V).$ Thus ${\rm Ind}_{H}^{G}({\tilde U}) = V$ is simple.

\medskip
Now ${\rm Res}^{G}_{H}(V) = {\tilde U} \oplus T$ for some projective $kH$-module $T,$ as $H$ is strongly
$p$-embedded in $G.$ Thus ${\rm Res}^{G}_{H}( V \otimes V^{*}) = k \oplus W$ for some projective
$kH$-module $W.$ By Green correspondence, $V \otimes V^{*} = k \oplus S$ for some projective $kG$-module 
$S.$ Thus $V$ is endotrivial.

\medskip
For the purpose of eliminating possible components $L$ of finite groups $G$ which contain
an elementary Abelian subgroup of order $p^{2},$ yet have a faithful endotrivial $kG$-module which is not
monomial, the following may be useful:

\medskip
\noindent {\bf Lemma 6:} \emph{ Let $L$ be a finite quasi-simple group of order divisible by $p$ which 
is a component of a finite group $G$ containing an elementary Abelian $p$-subgroup of order $p^{2}$
and having a faithful simple endotrivial $kG$-module, say $V.$ Let $Q$ be a Sylow $p$-subgroup of $G$
and $P = Q \cap L.$ Then there is faithful simple $kL$-module $U$ whose dimension is congruent to $\pm 1$ (mod $|P|$)
if $p$ is odd, or congruent to $\pm 1$ mod($\frac{|P|}{2}$) if $p =2.$ If $P$ is cyclic,  then $|P| > p$ and there is a 
faithful simple $kL$-module $U$ with ${\rm dim}_{k}(U) \equiv \pm 1$ (mod $p|P|).$}

\medskip
\noindent {\bf PROOF:} Notice that the hypotheses imply that $Z(L)$ is a $p^{\prime}$-group. In particular, notice that 
(by the Brauer-Suzuki theorem), $P$ is not generalized quaternion. We set $U = {\rm Res}^{G}_{L}(V),$ which is simple by Lemma 2 (and is endotrivial). Since ${\rm dim}_{k}(V)^{2} \equiv 1$ (mod $|Q|$), the stated congruences hold in all cases (notice that $P \neq Q$ by hypothesis when $P$ is cyclic). Also, $L$ contains an elementary
Abelian subgroup of order $p^{2}$ unless $p$ is odd and $P$ is cyclic. In that exceptional case, ${\bar L} = L/Z(L)$ admits an outer
automorphism $\alpha$ of order $p,$ since $F^{*}(G) = LZ(G)$ by Lemma 2 (and since $O_{p}(G) = 1)$. By the classification
of finite simple groups, $|P|>p$ in that case- we outline one argument: by a Frattini argument, we may suppose first that
$\alpha $ normalizes $N_{{\bar L}}({\bar P}).$ Another Frattini-type argument (using a Hall $p^{\prime}$-subgroup of
$N_{{\bar L}}({\bar P})$) allows us to suppose that $\alpha$ normalizes
a Hall $p^{\prime}$-subgroup ${\bar H}$ of $N_{{\bar L}}({\bar P}).$ If $|P| = p,$ then $\alpha$ is central in
$\langle \alpha \rangle N_{{\bar L}}({\bar P}\langle \alpha \rangle),$ which is a Sylow $p$-normalizer of ${\bar L}\langle \alpha \rangle.$
Hence $C_{\langle \alpha \rangle {\bar L}}(\alpha)$ controls strong $p$-fusion in $\langle \alpha \rangle {\bar L}.$
Then $[{\bar L},\alpha] = 1,$ since $O_{p^{\prime}}({\bar L}\langle \alpha \rangle) = 1$ and the odd analogue of
Glauberman's $Z^{*}$-theorem holds. This is a contradiction, as $\alpha$ induces a non-trivial outer automorphism
of ${\bar L}.$ Hence $|P| >p$ and the proof is complete.

\medskip
\paragraph{Acknowledgements:} We are grateful to C.W. Parker for helpful exchanges and to M. Geck for
computational assistance.

\end{document}